\documentclass[12pt, twoside]{amsart}
\usepackage[english]{babel}
\theoremstyle{plain}
\addtocounter{page}{1}

\def \R {{\mathbb {R}}}
\def \N {{\mathbb N}}
\def \C {{\mathbb C}}

\renewcommand{\a}{*}
\newcommand{\B}{\mathbf B}

\renewcommand{\d}{{\tt{d}}}

\renewcommand{\l}{\lambda}

\renewcommand{\r}{\rho}

\renewcommand{\Re}{\operatorname{Re}}
\renewcommand{\bar}{\overline}

\newcommand{\Ctr}{\mathcal C}

\newcommand{\Ctre}{\mathcal C_\epsilon}

\newcommand{\Ctrse}{\mathcal C^\sharp_\epsilon}

\newcommand{\Etre}{\mathcal E_\epsilon}

\newcommand{\h}{\mathfrak h}
\newcommand{\ha}{\mathfrak h^*}

\newcommand{\hp}{\mathcal{H}^p}

\renewcommand{\iota}{s}

\newcommand{\Rtre}{\mathcal R_\epsilon}

\newcommand{\Rtrse}{\mathcal R^\sharp_\epsilon}

\renewcommand{\S}{\mathcal S}

\newcommand{\Sl}{\mathcal S_{\!\,\lambda}}
\newcommand{\Sp}{\mathcal S_\varphi}

\newcommand{\Cn}{\mathbb C^n}

\newcommand{\dee}{\partial}

\newcommand{\deebar}{\overline\partial}
\newcommand{\de}{\mathrm{d}}

\renewcommand{\N}{\mathcal N}
\newcommand{\bndry}{b}

\numberwithin{equation}{section}

\newcommand{\Aop}{\mathcal A^{(\sigma)}}
\newcommand{\Aope}{\mathcal A^{\!(\lambda)}_\epsilon}

\newcommand{\Sop}{\mathcal S}

\renewcommand{\dee}{\partial}
\renewcommand{\deebar}{\overline\partial}
\renewcommand{\de}{d}

\title[Harmonic Analysis techniques in SCV]{Harmonic Analysis Techniques in Several Complex Variables
\\
Su un' applicazione dell' Analisi Armonica reale all'analisi complessa in piu' variabili
}
\author{Loredana Lanzani$^*$}
\thanks{$^*$ Supported in part by the National Science Foundation, award DMS-1503612}
\address{Department of Mathematics\\Syracuse University\\215 Carnegie Bldg.\\Syracuse, NY 13244-1150 USA
}
\email{llanzani@syr.edu
}
\thanks{Bruno Pini Mathematical Analysis Seminar, Vol. 1 (2014) pp.   
\\
Dipartimento di Matematica, Universit\`{a} di Bologna
\\
ISSN 2240-2829}

\begin{document}

\begin{abstract}
{We give a survey of recent joint work with E. M. Stein (Princeton University) concerning the application of 
suitable versions of the T(1)-theorem technique to the study of orthogonal projections onto
the Hardy and Bergman spaces of holomorphic functions for domains
  with minimal boundary regularity.
}

\medskip \noindent {\sc{Sunto.}} 
{Questo resoconto offre una sintesi di una serie di recenti collaborazioni con E. M Stein (Princeton University)
sull'applicazione del celeberrimo teorema T(1) allo studio delle proiezioni ortogonali sugli spazi 
di Hardy e di Bergman 
per funzioni olomorfe su domini
dotati di minima regolarita' al bordo.
}


 \thanks{THIS IS AN ELECTRONIC REPRINT OF THE ORIGINAL ARTICLE THAT APPEARED IN THE {\em BRUNO PINI MATH. ANALYSIS SEMINAR} U. OF BOLOGNA (2014), 83-110 (ISSN: ISSN 2240-2829). THIS REPRINT DIFFERS FROM THE ORIGINAL IN PAGINATION AND TYPOGRAPHICAL DETAIL}

\end{abstract}
\maketitle
\section{
Introduction}
This is a summary of recent work 
\cite{LS-1}-\cite{LS-6}
 concerning the
 $L^p$-regularity properties of {\em orthogonal} projections (Bergman projection, Szeg\H o projection) onto $L^2$-closed subspaces of {\em holomorphic functions} (Bergman space, holomorphic Hardy space) for bounded domains $D\subset \Cn$ with {\em minimal} boundary regularity. 
 Regularity properties of the Szeg\H o 
   and Bergman projections, 
   in particular $L^p$-regularity, have been the object of considerable interest for more that 40 years. 
When the boundary of the domain $D$ is sufficiently smooth, decisive results were obtained in the following settings:
{\em (a)}, when $D$ is strongly pseudoconvex \cite{KS-1}, \cite{L};
{\em (b)}, when $D \subset \mathbb C^2$ and its boundary is of finite type
 \cite{Mc-1}, \cite{NRSW};
{\em (c)}, when $D \subset \mathbb C^n$ is convex and its boundary is of finite type
 \cite{Mc-2}, \cite{MS-2}; and
 {\em (d)}, when 
 $D \subset \mathbb C^n$ is of finite type and its Levi form is diagonalizable
 \cite{CD}. 
 Related results include
  \cite{AS-1}, \cite{Ba}, \cite{BaSa}, \cite{BaVa}, 
  \cite{BoLo}, 
   \cite{EL},
   \cite{FH-1}-\cite{FH-5},
  \cite {Ha}, \cite{Hans}, \cite{HNW}, 
 \cite{KP}, \cite{NRSW}, \cite{Ro-1}, \cite{Ro-2}, \cite{St}.
 
 It should be noted that several among these works depend on good estimates or explicit formulas for the Szeg\H o or Bergman kernels. 
 In our non-smooth setting these are unavailable and we have to proceed via a different framework, by
pursuing   a theory of  singular integral operators with {\em holomorphic} kernel
   that blends the complex structure of the ambient domain with the
     Calder\`on-Zygmund
     theory for singular integrals on non-smooth domains in $\mathbb R^{2n}$. 
Our present
task is
 to highlight the main threads linking the various themes in
  \cite{LS-1}-\cite{LS-6} and convey a general idea of the methods of proof (and at times we will 
 sacrifice technical detail in favor of a more streamlined exposition). While most of the proofs are deferred to \cite{LS-1}-\cite{LS-6}, here we indicate references to the specific statements therein. 
 \medskip
 
 {\bf Aknowledgment.} I am grateful to the organizers and participants of the {\em Bruno Pini Mathematical Analysis Seminar} for the kind hospitality and lively discussions.

\section{The Szeg\H o projection}
\subsection{Motivation and context} Our starting point is the seminal work by Calder\`on \cite{C}, Coifman-McIntosh and Meyer \cite{CMM} and David \cite{D}  on the $L^p(\Gamma)$-regularity of the classical Cauchy integral for a planar curve $\Gamma\subset \C$, in the situation when $\Gamma$ is the boundary of a 
domain $D\subset \C$ (and we will write $\Gamma=\bndry D$):
\begin{equation}\label{E:Ctr}
\Ctr f(z) = \frac{1}{2\pi i}\!\int\limits_{w\in\bndry D}\!\!\!\!\! f(w)\, \frac{dw}{w-z} ,\quad z \in\bar D.
\end{equation}
 For $z\in \bndry D$ we interpret \eqref{E:Ctr} as a singular integral in the ``principal value'' sense, see
 \cite[(1.1)]{C-2}.
The situation when $\bndry D$ is of class $C^{1, \alpha}$ (with $\alpha>0$) can be easily reduced to the classical setting of the Hilbert transform operator \cite[Section 1.1, Example 8]{C-2}. However dealing with the case when $\bndry D$ is of class $C^1$ and more generally Lipschitz, required new ideas that ultimately led to the so-called ``T(1)-theorem'' technique \cite{DJ} for a more general class of singular integral operators\footnote{that include the
  the Cauchy integral \eqref{E:Ctr} as a prototype.}, and to applications to the study of analytic capacity \cite{To-1} as well as the solution of the Vitushkin conjecture \cite{To-2}.  
 In the setting of higher dimension (that is, 
for a Lipschitz domain $D\subset \R^N$ with $N\geq 2$), the Cauchy integral and the
 related singular integral operators collectively known as {\em boundary layer potentials}\footnote{namely, the single layer potential and the double layer potential operators, see \cite[Section 2.2]{Ke} and references therein.} provide the solution to various boundary value problems for harmonic functions.
Here we are especially interested in
the {\em $L^p$-Dirichlet
 problem for harmonic functions}:
given $u: \bndry D \to \R$ with $u\in L^p(\bndry D, d\sigma)$,  where $d\sigma$ is the induced Lebesgue measure on $\bndry D$, 
find\ $U: D\to\R$\ such that
\begin{equation} 
\label{E:BVP-harm}
\left\{
\begin{array}{rll}
U_{xx}+ U_{yy}\ =:\ \Delta U (z)=&\!\!\!0,
  &\text{if}\ z \in D\ \\
\lim\limits_{z\to w} U(z) =:U^+(w)=&\!\!\!u(w), &\text{if}\  w\in \bndry D\\
 \|\N(U)\|_{L^p(\bndry D, d\sigma)}\leq &C \|u\|_{L^p(\bndry D, d\sigma)}\, ,&
 \end{array}
 \right.
\end{equation}
where the limit that defines the boundary value $U^+$ is to be suitably interpreted (for instance, as a 
``non-tangential limit'' \cite[page 24]{Ke}) and $\N(U)$ denotes the so-called ``non-tangential maximal function'' for $U$, see \cite[page 13]{Ke} and references therein. The solution of 
\eqref{E:BVP-harm}
can be expressed in terms of the aforementioned boundary layer potentials acting on the data $u$. As it turned out, 
 the size of the $p$-range for which existence, uniqueness and $L^p$-regularity of the solution occur, is related to the size of the Lipschitz constant of the domain \cite[Theorem 2.2.2]{Ke}. (We may think of $C^1$- or smoother domains as having Lipschitz constant equal to zero.)

\subsection{Regularity of the Szeg\H o projection: statement of the problem}

We are interested in 
the {\em holomorphic} analog of the 
 problem  \eqref{E:BVP-harm} for domains with minimal regularity, which we presently recall in the situation when $D\subset \C\equiv \R^2$ (the planar setting). The $L^p$-Dirichlet problem for {\em holomorphic} functions on a planar domain $D\subset \C$ is stated as follows.
Given $g: \bndry D \to \C$ in $L^p(\bndry D, d\sigma)$,
find $G: D\to \C$ such that
\begin{equation} 
\label{E:BVP-hol}
\left\{
\begin{array}{rll}
1/2(G_x-i G_y)\ =:\ \deebar G (z)=&\!\!\!0,
  &\text{if}\ z \in D\\
\!\!\!\!  \!\!\!\! \lim\limits_{z\to w} G(z)=: G^+(w)= &\!\!\! g(w), &\text{if}\  w\in \bndry D\\
\|\N(G)\|_{L^p(\bndry D, d\sigma)}\leq &C \|g\|_{L^p(\bndry D, d\sigma)}, &
 \end{array}
 \right .
\end{equation}
where we adopt the convention that $z\in D$ is expressed as $z=x+iy$. It is clear that if $G$ solves \eqref{E:BVP-hol} with data $g=u+iv$ then e.g.,  $U:= \mbox{Re}\, G$ solves \eqref{E:BVP-harm} with data $u$.
  However,
  in contrast with the situation for
    \eqref{E:BVP-harm},
the natural data space for \eqref{E:BVP-hol} is not $L^p(\bndry D, d\sigma) +i L^p (\bndry D, \sigma)$:
it is instead the {Hardy Space of holomorphic functions} {$H^p(D)$}  (aka {\it Smirnov Class}) 
$$
{H^p(D)}:=\bigg\{F \, \bigg|\, \deebar F (z) =0, z\in D,\ \sup\limits_{\epsilon>0}
\!\!\int\limits_{z\in\bndry D_\epsilon}\!\!\!\!\!|F(z)|^p\de\sigma_\epsilon(z)
<+\infty\bigg\}
$$
where $\{D_\epsilon\}_\epsilon$ is any family of (say, rectifiable) subdomains of $D$ with $D_\epsilon \uparrow D$.

 In fact $H^p(D)$ can be identified with a {\em proper} subspace of $L^p(\bndry D, d\sigma) +i L^p (\bndry D, d\sigma)$ which we denote $\hp (\bndry D, d\sigma)$.
  More precisely, we invoke the well-known fact that functions in  $H^p(D)$ have non-tangential limits that belong to  $L^p(\bndry D, d\sigma)$, see \cite{Du} and \cite{S-1}, and then we identify $\hp(\bndry D, d\sigma)$ with the space $\{F^+\ |\ F\in H^p(D)\}$.
Returning to the Dirichlet problem for holomorphic functions,  one thus needs $g = F^+$ 
 for some $F\in H^p(D)$  and if this is the case then
$G:= F$ solves \eqref{E:BVP-hol}.

The holomorphic Hardy space $\hp(\bndry D, d\sigma)$ is a {closed} subspace
of $L^p(\bndry D, d\sigma)$ (a fact that can be seen e.g., by applying the Cauchy formula on small discs + the co-area formula \cite{Fe}); thus for the exponent
{$p=2$} the theory of Hilbert spaces grants the existence of
a unique, {\em orthogonal} projection $\Sop$: $L^2(\bndry D, d\sigma)\to H^2(\bndry D, d\sigma)$,
known as the {\em Szeg\H o projection},
which is
  a singular integral operator 
 characterized by the following 
properties: 
$$
\Sop^2=\Sop,\quad  
 \Sop^\a=\Sop, \quad
\|\Sop\|_{L^2 (\bndry D, \sigma)\to L^2(\bndry D, d\sigma)}=1\, .
$$
(Here $\Sop^\a$ denotes the adjoint of $\Sop$ taken with respect to the inner product
in $L^2(\bndry D, d\sigma)$.) 

These properties in particular indicate that the  Szeg\H o projection is the natural solution operator for 
\eqref{E:BVP-hol} in the case when $p=2$. 
On the other hand, solving \eqref{E:BVP-hol} in $L^p(\bndry D, d\sigma)$ for $p\neq 2$ is a much harder
problem, and one that 
is ultimately related to the
  \smallskip
    
     \underline{\bf $L^p$-Regularity problem for $\Sop$}: {\em under {\em minimal} assumptions on $D$, find the
 {\em maximal} exponent $P=P(D) \in [2, +\infty]$ such that
  $\Sop: 
L^p(\bndry D, d\sigma)\to L^p(\bndry D, d\sigma)$ is bounded {\em for all $P'< p< P$.}
}
\smallskip

By symmetry considerations (the fact that $\Sop^*=\Sop$) we have that $P$ and $P'$ must be conjugate exponents (namely, $1/P + 1/P'=1$).

We point out that the problem \eqref{E:BVP-hol} and
the $L^p$-regularity problem for $\Sop$ can also be stated in higher dimension, that is for $D\subset\Cn$: in this setting the quantity $\deebar G$ is interpreted as a differential form of type $(0, 1)$ and the condition that $\deebar G =0$ is then equivalent to  the requirement that$$
G_{x_j}-i G_{y_j}=0\quad \mathrm{for\ all}\quad j=1,\ldots, n
$$
where $x_j+iy_j = z_j$ with $j=1,\ldots n$ and $z=(z_1, \ldots, z_n)\in D$.  (Everything else in \eqref{E:BVP-hol}, also the definitions pertaining the Hardy spaces
$\hp(\bndry D, d\sigma)$ and the $L^p$-regularity problem for $\Sop$,
are meaningful regardless of the size of the dimension\footnote{A more general version of these problem can be stated in which the data $g$ is a differential form of degree $0\leq r\leq 2n-1$ and  includes \eqref{E:BVP-hol} as the special case: $r=0$, but we will not pursue it here.}).

\subsection{Regularity of the Szeg\H o projection: case of planar domains.}\label{SS:2.3}
 It turns out that in complex dimension $n=1$, that is for a bounded and simply connected domain $D\subset \C$, 
the size of the maximal interval $(P', P)$ is related to the boundary regularity of $D$. Specifically, we have the following results:

\begin{itemize}
\item[1.]  If $D\subset \mathbb C$ is {Vanishing Chord-Arc} (e.g., $D$ of {class $C^1$}), then
{$ P=+\infty $}, 
\cite[Theorem 2.1 (1)]{LS-1} (see also \cite{La-1}).
\item [2.] If $D\subset \mathbb C$ is {Lipschitz with constant $M$}, then
$${\displaystyle{P = 2\left(1+\frac{\pi}{2\arctan M}\right)>\ 4}}\, ,$$ 
and the interval determined by such $P$ is optimal within the Lipschitz category, 
\cite[Theorem 2.1 (2), and page 69]{LS-1}.
\item [3.] If $D\subset \mathbb C$ is a {rectifiable local graph}, then
{$ P =4 $}, \cite[Theorem 2.1 (3)]{LS-1}.
\item [4.] If $D\subset \mathbb C$ is {Ahlfors-Regular}, then {$P=2+\epsilon$ for some $\epsilon =\epsilon (D)>0$}, \cite[Theorem 2.1 (4)]{LS-1}.
\item [5.] There is a rectifiable domain $D_0\subset \mathbb C$, see  \cite{Be}, such that 
$$\Sop: L^p(\bndry D_0, \sigma)\to L^p(\bndry D_0, \sigma) \iff {p=2}.$$
\end{itemize}

The methods of proof for all these results
 rest on the existence of a conformal map $\psi: \mathbb D_1(0)=\{|z|<1\}\to D$
 (namely, the original problem for $\Sop$ is reduced to a {\em weighted} problem for $\S_0$= the Szeg\H o projection of $\mathbb D_1(0)$ with weight $\omega = |\psi'|^{1-p/2}$ to which one may apply the theory of Muckenhoupt  \cite{S-2})
  and thus
  are not applicable to
     higher dimension that is, to the situation when $D\subset \Cn$ and $n\geq 2$. 
     
     On the other hand, 
   item 1. can also be studied via
     a conformal map-free argument that relies on the
     comparison of $\S$ with the Cauchy integral $\Ctr$. We point out that for the Cauchy integral boundedness in $L^2$ {\em implies} boundedness in $L^p$
  for $1<p<\infty$, see \cite{S-2}, 
and so in general we have $\Sop\neq \Ctr$.

The approach to the analysis of the Szeg\H o projection that we are about to describe was first formulated for the case
when $D\subset \C$ is smooth, see \cite{KS-1} and \cite{KS-2}, and
the comparison of $\Ctr$ and $\Sop$ 
hinged on the following facts:
 \begin{itemize}
\item[\tt (a.)]  Each of $\Ctr$ and $\Sop$ is a {\em projection\footnote{ that is, $\Sop\circ\Sop=\Sop$ and $\Ctr\circ\Ctr=\Ctr$.}}: $L^2(\bndry D, d\sigma)\to \mathcal H^2 (\bndry D, d\sigma)$.
\item [\tt (b.)] $\Sop$ is self-adjoint while $\Ctr$ (in general) is not\footnote{unless $D$ is a disc, see \cite{KS-2}.}.
\item[\tt (c.)] Each of $\Ctr$ and $\Sop$ is {\em bounded}: $L^2(\bndry D, d\sigma)\to L^2(\bndry D, d\sigma)$.
\item [\tt(d.)] The kernel of the operator $\Aop := \Ctr^* -\Ctr$, where $\Ctr^*$ denotes the formal adjoint of $\Ctr$ in $L^2(\bndry D, d\sigma)$, is {\em ``small''} if $D$ is sufficiently smooth (a cancellation of singularities occurs by performing a second-order Taylor expansion at $w=z$).
  \end{itemize}
  Then one has the following identities on $L^2(\bndry D, d\sigma)$:
  $$
  \Sop\Ctr = \Ctr\quad \mathrm{and}\quad  \Ctr\Sop=\Sop\, ,\quad \mathrm{by\ item\ \tt(a.)}\,.
  $$
  Taking $L^2(\bndry D, d\sigma)$-adjoints of the second identity above, we get
  \begin{equation*}
  \Sop\Ctr^*=\Sop\, ,\quad\mathrm{by \ item\ \tt(b.).}
  \end{equation*}
  Subtracting the first of the two identities above from the latter we obtain
  \begin{equation}\label{E:KS-0}
  \Sop\,[I-\Aop] = \Ctr \quad \mbox{in}\quad L^2(\bndry D, d\sigma)
    \end{equation}
  where $I$ denotes the identity operator: $L^2(\bndry D, d\sigma)\to L^2(\bndry D, d\sigma)$.
   Now using {\tt(c.)} and the fact that
    $(\Aop)^* = - \Aop$ (recall that $\Aop = \Ctr^*-\Ctr$)
  it is not hard to see that the operator 
  $I-\Aop$ is invertible on $L^2(\bndry D, d\sigma)$ with bounded inverse,
     and we conclude that the identity
  \begin{equation}\label{E:KS-1}
  \Sop = \Ctr\, [I-\Aop]^{-1}\quad \mathrm{holds\ in}\quad L^2(\bndry D, d\sigma).
  \end{equation}
  However, by item {\tt (d.)}
  (which holds if $D$ is smooth)
  the operator $\Aop$ is in fact compact in $L^p(\bndry D, d\sigma)$ for $1<p<\infty$, and by the closed graph theorem it follows that $I-\Aop$ is invertible in $L^p(\bndry D, d\sigma)$ with bounded inverse, see \cite[page 65]{LS-1}.  It follows that the right-hand side of
    \eqref{E:KS-1} is a well-defined and bounded operator: $L^p(\bndry D, d\sigma)\to L^p(\bndry D, d\sigma)$ for $1<p<\infty$, 
     and we conclude from the above
  that $\Sop$ extends to a bounded operator on $L^p(\bndry D, d\sigma)$ for $1<p<\infty$, thus solving the $L^p$-regularity problem for $\Sop$ with $P=\infty$,  whenever $D\subset\C$ is smooth.
  
  We remark that the steps {\tt(a.)} --  {\tt(d.)} can be stated for {\em any} positive boundary measure $\mu$ (not just the induced Lebesgue measure $\sigma$)
  provided the orthogonal projection $\S \equiv \S_\mu$ is defined with respect to the duality
   induced by the measure $\mu$, namely
  $$
  (f, g)=\int\limits_{w\in\bndry D}\!\! f(w)\bar g(w)\, d\mu (w).
  $$
  
  \subsection{Regularity of the Szeg\H o projection: dimension-induced obstructions.}
  \label{SS:2.4}
The procedure described in the previous section is, in principle, dimension-free in the sense that it relies on the existence of ``some'' operator $\Ctr$ that satisfies the four conditions {\tt (a.)} through {\tt (d.)}.
 In the setting of Section \ref{SS:2.3} (that is when $D\subset \C$ and $D$ is sufficiently smooth) one  takes $\Ctr$ to be the Cauchy integral
 \eqref{E:Ctr}, and
the proof of the crucial item {\tt (a.)} then
  rests on the following two 
  features of $\Ctr$:
   
  \begin{itemize}
  \item[\tt(i.)] The fact that 
     {\em Cauchy kernel} $C(w, z)$ 
  (that is the kernel of $\Ctr$)
   is {\em universal} in the sense that its dependence on the domain $D$ is effected only through the inclusion
   
  \centerline{$j:\bndry D \hookrightarrow \C\, . $}
   
 \noindent Specifically, we have
  \begin{equation}\label{E:Cauchy-ker-1}
  C(w, z) = \frac{1}{2\pi i}\,j^*\!\left(\frac{\de w}{w-z}\right) \quad w, z\in \C\times \C, 
\ \  w\neq z\, ,
  \end{equation}
  where $j^*$ is the so-called {\em pull-back by the inclusion map}, see e.g.,  \cite[Section III.1.5]{Ra-1}.
  \item[]
  \item[\tt(ii.)] The fact 
  that the {\em Cauchy kernel function} $1/(w-z)$ is (obviously) {\em holomorphic} in the parameter $z\in D$ whenever $w\in \C\setminus \bar{D}$, in particular for each fixed $w\in\bndry D$.
 \end{itemize}
 In higher dimension both of these properties become highly problematic
 as
  the only known {\em universal\,}
   reproducing kernel
  is the {\em Bochner-Martinelli} kernel:
  \begin{equation}\label{E:BM}
H(w, z) = 
 \frac{(n-1)!}{(2\pi i)^n}\sum\limits_{j=1}^n j^*\!\left(\frac{\bar w_j -\bar z_j}{|w-z|^{2n}}\,
\de w_j\bigwedge\limits_{\nu\neq j}\de\bar w_\nu \wedge \de w_\nu\right)\!,\ \ w, z\in\Cn\times\Cn\,, w\neq z,
\end{equation}
see e.g., 
  \cite[Lemma IV.1.5 (a)]{Ra-1}. It is clear that $H(w, z)=C(w, z)$ when $n=1$, because in such case the coefficient in $H(w, z)$ is just $(\bar w -\bar z)/|w-z|^2 = 1/(w-z)$. 

On the other hand, when $n\geq 2$ the coefficients of the
kernel 
\eqref{E:BM}
are obviously nowhere holomorphic, thus $H(w, z)$ is of no use in the 
analysis of
 the Szeg\H o projection described in the previous section\footnote{The Bochner Martinelli integral for a general domain $D$ cannot satisfy item {\tt (a.)}.}: there is no canonical, higher dimensional {\em holomorphic} analog of the Cauchy kernel \eqref{E:Cauchy-ker-1}. Instead,  one has to look into ad-hoc constructions that are tailored to certain specific features of the domain. More precisely,  the existence of a higher-dimensional holomorphic analog of $C(w, z)$
 is intimately connected with a geometric constraint on the domain, namely the requirement that $D$ be {\em pseudoconvex} \cite[Section II.2.10]{Ra-1} or, equivalently, that $D$ be
 a so-called {\em weak (or local) domain of holomorphy} \cite[Section II.2.1]{Ra-1}: for any $w\in \bndry D$ there must be a function $f_w(z)$ that is holomorphic in $z\in D$ but cannot be extended holomorphically past $w$.
  While any planar domain $D\subset \C$ is obviously a weak domain of holomorphy\footnote{for any $w\in\bndry D$,  simply take  $f_w(z)$ to be the Cauchy kernel function, i.e. $f_w(z) =1/(w-z)$, $z\in D$.}, 
  there are domains $D\subset \Cn$ ($n\geq 2$) with the property that {\em any}  function  holomorphic in $D$ can be holomorphically extended to a larger domain $\Omega\supset D$ \cite[Lemma II.2.2]{Ra-1}.

\subsection{ Regularity of the Szeg\H o projection in higher dimension: the case of smooth domains}\label{SS:2.5}
  The Cauchy-Fantappi\`e theory (see \cite[Section 4]{LS-3} and references therein) provides an algebraic
  framework
  to 
    construct
   explicit, higher dimensional holomorphic
  analogues of the Cauchy kernel \eqref{E:Cauchy-ker-1} 
    for any bounded, {\em strongly} (equivalently, {\em strictly}) pseudoconvex domain $D\subset \Cn$, see \cite[Section II.2.8]{Ra-1}. The kernel construction and the proof of the 
       corresponding conditions 
     {\tt (a.)} through {\tt(d.)} were first carried out in \cite{Hen-1}, \cite{KS-1} and \cite{R}
      and dealt with the case when the strongly pseudoconvex domain $D$  is smooth (of class $C^3$ or better). In this section we describe the construction in such setting (see also \cite{Ker} and \cite{LS-3}).
  
  For $D$ strongly pseudoconvex we write $D=\{\rho (z)<0\}$ where $\rho:\Cn \to \R$ is a strictly plurisubharmonic {\em defining function} for $D$ (see \cite[Sections II.2.3 and  II.2.7 ]{Ra-1}), which is taken to be of class $C^3$ or better. For fixed $z\in D$, we consider the following differential form of type $(1, 0)$ in the variable $w$
  $$
\eta(w, z) =   \sum\limits_{j=1}^n\eta_j(w, z)\, dw_j
  $$
  where we have set
  \begin{equation}\label{E:eta}
\eta_j(w, z) = \ \chi_0(w, z)
{\left(\frac{\dee\rho}{\dee\zeta_j}(w) -
\frac12\sum\limits_{i=1}^n
\frac{\dee^2\rho (w)}
{\dee\zeta_i\dee\zeta_j}
(w_i-z_i)\right)} + (1-\chi_0(w, z))(\bar w_j-\bar z_j)
\end{equation}
  with $\chi_0$ a smooth cutoff function supported in $\{|w-z|<\delta\}$.

\noindent    Now $\eta (w, z)$ is a {\em generating form at $w$} in the sense that the complex-valued function of $z$
$$
\langle\eta (w, z), w-z\rangle:=\sum\limits_{j=1}^n\eta_j(w, z)\,(w_j-z_j)
$$
is bounded below by $|w-z|^2$  for any $z\in\bar D$, see \cite[Section 4]{LS-3}. More precisely we have
\begin{equation}\label{E:est-1}
\Re \langle\eta (w, z), w-z\rangle\geq c |w -z|^2, \quad w\in \bndry D,\ z\in \bar D\, .
\end{equation}
 (We point out that the validity of this inequality when $z$ is close to $w$ is a consequence of the strict plurisubharmonicity of $\rho$, see \cite[Lemma 4]{LS-3}.) The Cauchy-Fantappi\`e theory then grants that the 
 kernel
 \begin{equation}\label{E:C-til}
 \widetilde C(w, z)= \frac{1}{(2\pi i)^n}
{\frac{\eta\wedge(\deebar_w\eta)^{n-1}(w, z)}
{\langle\eta(w, z), w-z\rangle^n}}
 \end{equation}
  reproduces holomorphic functions (more precisely, the induced singular integral fixes the space $\mathcal H^2(\bndry D, d\sigma)$) and it is clear that $\widetilde C(w, z)$   satisfies property {\tt (b.)}, see \cite[Section 4]{LS-3} and references therein\footnote{Roughly speaking, one wants $C(w, z) \neq \bar{C(z, w)}$, which is the case whenever $D\neq \{|z|^2<1\}$.}. On the other hand, it is apparent from \eqref{E:eta} that, as a function of $z\in D$, this kernel is holomorphic only for $z$ near $w\in\bndry D$. Thus, in order to achieve the crucial condition {\tt (a.)} one needs to modify $\widetilde C(w, z)$ by adding a correction term that will make the kernel globally holomorphic: 
 \begin{equation}\label{E:HCK-3}
  C(w, z)\ =\
  \widetilde C(w, z)
\ + \ C_\rho (w, z)\, .
  \end{equation}
  The correction $C_\rho (w, z)$ is obtained either by solving a $\deebar$-problem (in the $z$-variable)
   on a strongly pseudoconvex, smooth domain $\Omega$ that contains $\bar D$, see \cite{KS-1} and \cite[Section 8]{LS-3}, or by solving a Cousin problem as in  \cite{Hen-1} and \cite{R}. The resulting kernel \eqref{E:HCK-3} will 
  be globally holomorphic and the corresponding operator (still denoted $\Ctr$)  will 
    satisfy properties {\tt (a.)} and {\tt (b.)}.   
 
  We point out that the procedure described up to this point 
   can be carried out
    under the weaker assumption that
  the domain (that is the defining function $\rho$) be of class $C^2$ \cite[Section V.1.1]{Ra-1}: 
  it was in order  to prove the remaining properties {\tt (c.)} and {\tt (d.)} 
   that one needed
  more regularity.
     Specifically, in the setting of \cite{KS-1} one needed to assume that  $\rho$ be of class $C^3$, and the proof of item {\tt (c.)}
   (the $L^2(\bndry D, d\sigma)$-regularity of the holomorphic Cauchy integral $\Ctr$)
          could then be achieved  via an ``osculation by model domain'' technique. The basic idea is that
     there is a strongly pseudoconvex and smooth ``model'' domain $D_0$ for which
      the operator $\Ctr_0$ constructed as in \eqref{E:C-til} and \eqref{E:HCK-3}
             takes an especially simple form, and the validity of property {\tt (c.)}
    for such a $\Ctr_0$ is easily verified by direct inspection\footnote{In fact $D_0$ is the {\em Siegel upper half space}: $\{z=(z', x_n+iy_n)\in \Cn\ |\  |y_n|>|z'|^2\}$, where $z'=(z_1,\ldots, z_{n-1})$.}.  On the other hand, if $D$ is strongly pseudoconvex and of class $C^3$ then at any boundary point it is osculated by a copy of $D_0$ with small error. Furthermore, one may write the operator $\Ctr$ (for the original domain $D$) as the sum of $\Ctr_0$ (the corresponding operator for the model domain $D_0$) plus  the ``error'' operator $\Ctr -\Ctr_0$, and if $D$ is of class $C^3$ the error operator is easily seen to be bounded, thus concluding the proof of {\tt (c.)}. Finally, a 2nd-order Taylor expansion of $\rho$ in the variable $z$ about the point $w$ shows that the kernel of $\Aop = \Ctr^*-\Ctr$ has the ``smallness'' property {\tt (d.)} whenever $D$ is of class $C^3$ (which allows for good control on the tail of the expansion). 
    
 Having constructed an operator  $\Ctr$ that satisfies the four properties {\tt (a.)} through {\tt (d.)}, one proceeds as in Section \ref{SS:2.3} to conclude  that the $L^p(\bndry D, d\sigma)$-regularity problem for $\Sop$ is  solved with $P=\infty$,
 whenever $D$ is strongly pseudoconvex and of class $C^3$, see \cite{Ker} and  \cite{KS-1}.
     
     We point out that the methods of proof for items {\tt (c.)} and {\tt (d.)} as described above break down as soon as the regularity of  $D$  is below the class $C^3$. (The ``error'' operator $\Ctr-\Ctr_0$ that occurred the proof of {\tt (c.)} can no longer be controlled, whereas for {\tt (d.)} there is no control on the size of the tail in the aforementioned Taylor expansion.)
     
 \subsection{Higher dimensional holomorphic kernels for non-smooth domains: kernel construction}\label{SS:2.6} We now describe the results  in
\cite{LS-5}. As we have seen, a natural requirement for the existence of a holomorphic Cauchy-type kernel \eqref{E:HCK-3} is that the domain be strongly pseudoconvex, which is a condition that  essentially involves two degrees of differentiability of the boundary of the domain. As a result, the threshold of smoothness for a strongly pseudoconvex domain should be the class $C^2$ (as opposed to the class $C^1$ for a planar domain):
as before, we take $\rho$ to be  a strictly plurisubharmonic defining function for $D$, however now $\rho$ is merely of class $C^2$. 
 To make up for the lack of differentiability of those
 second derivatives of $\rho$ that occurred in the definition  of the generating form $\eta$, see \eqref{E:eta}, we ``borrow some regularity'' by considering families of  functions $\{\tau^{(\epsilon)}_{i,j}\}_\epsilon$
of class $C^2$ such that
\begin{equation}\label{E:est-2}
\sup\limits_{w\in \bndry D}\left |\frac{\dee^2\rho (w)}
{\dee\zeta_i\dee\zeta_j} - \tau^{(\epsilon)}_{i,j}(w)\right |<\epsilon\, ,\quad i, j=1,\ldots, n
\end{equation}
for any $0<\epsilon \leq \epsilon_0$, where the size of $\epsilon_0$ is determined by the
the strict plurisubharmonicity of $\rho$, see \cite[(2.26)]{Ra-1}.

One then sets  $$
\eta^{(\epsilon)}(w, z) =   \sum\limits_{j=1}^n\eta^{(\epsilon)}_j(w, z)\, dw_j
  $$
with
  \begin{equation}\label{E:eta-eps}
\eta^{(\epsilon)}_j(w, z) = \ \chi_0(w, z)
{\left(\frac{\dee\rho}{\dee\zeta_j}(w) -
\tau^{(\epsilon)}_{i,j}(w)
(w_i-z_i)\right)} + (1-\chi_0(w, z))(\bar w_j-\bar z_j)
\end{equation}
  where $\chi_0$ is again a smooth cutoff function supported in $\{|w-z|<\delta\}$.
  
  It follows from \eqref{E:est-2} and \eqref{E:est-1} that

\begin{equation}\label{E:est-3}
\Re \langle\eta^{(\epsilon)} (w, z), w-z\rangle\geq c |w -z|^2
\end{equation}
and also
\begin{equation}\label{E:est-4}
\left |\langle\eta^{(\epsilon)} (w, z), w-z\rangle\right |\approx 
\left |\langle\eta (w, z), w-z\rangle\right |
\end{equation}
whenever $z\in\bar D$ and $w\in \bndry D$, and for any $0<\epsilon\leq \epsilon_0$,
 with the constant $c$ in \eqref{E:est-3} and the implied constants in \eqref{E:est-4}  independent of $\epsilon$, see \cite[Part I]{LS-5}.
It follows from \eqref{E:est-3} that, as was the case for $\eta$ in the previous section, each $\eta^{(\epsilon)}$ is a  generating form (for any $0<\epsilon<\epsilon_0$). Thus, the resulting Cauchy-Fantappi\` e kernels
\begin{equation}\label{LS-kere}
 \widetilde C^{(\epsilon)}(w, z)= \frac{1}{(2\pi i)^n}
{\frac{\eta^{(\epsilon)}\wedge(\deebar_w\eta^{(\epsilon)})^{n-1}(w, z)}
{\langle\eta^{(\epsilon)}(w, z), w-z\rangle^n}}
 \end{equation}
 have the reproducing property {\tt (a.)} but as before, are only locally holomorphic (for $z\in D$ near $w\in\bndry D$). To achieve global holomorphicity one again has to solve a $\deebar$-problem to produce suitable correction terms:
 \begin{equation}\label{E:LS-ker}
 C^{(\epsilon)}(w, z)= \frac{1}{(2\pi i)^n}
{\frac{\eta^{(\epsilon)}\wedge(\deebar_w\eta^{(\epsilon)})^{n-1}(w, z)}
{\langle\eta^{(\epsilon)}(w, z), w-z\rangle^n}} + C_\rho^{(\epsilon)}(w, z).
 \end{equation}
What matters here is that each of the corrections $C_\rho^{(\epsilon)}(w, z)$ satisfies a 
uniform bound which is independent of $\epsilon$, namely
$$
\sup\limits_{(w, z)\in \bndry D\times \bar D}|C_\rho^{(\epsilon)}(w, z)|\leq C\quad \mbox{for \ any}\quad 0<\epsilon\leq \epsilon_0.
$$
We let $\{\Ctre\}_\epsilon$ denote the resulting family of (globally) holomorphic Cauchy-type integral operators. It is clear from the above that each $\Ctre$ satisfies conditions {\tt (a.)} and {\tt (b.)} in Section \ref{SS:2.3}, for any $0<\epsilon<\epsilon_0$. 
 \subsection{$L^2(\bndry D)$-regularity of the $\Ctre$'s: preliminary observations}\label{SS:2.7}
  It turns out that the ``borrowed regularity'' \eqref{E:est-2} is not good enough to 
 prove $L^2$-boundedness of the $\Ctre$'s by the ``osculation by model domain'' method that was described in Section \ref{SS:2.5} (there is a problem with controlling the error 
 $\Ctre -\Ctr_0$, so 
  regularity for the $\Ctre$'s cannot be 
 deduced from the corresponding result for $\Ctr_0$) and we must proceed by a different route, namely, by the ``$T(1)$-theorem technique''. To this end, we make a number of preliminary observations.
 
$\bullet$ Our first observation \cite[Part I]{LS-5} is that there is an ad-hoc measure
  for $\bndry D$, which we call {\em the Leray-Levi measure $\l$}, that is better suited to study
  the $\Ctre$'s than the induced Lebesgue measure $\sigma$.

  More precisely, we set
  \begin{equation*}
d\l (w) \!=\! (2\pi i)^{-n}\! 
j^*\!\left(\dee\rho\!\wedge\!(\deebar\dee\r)^{n-1}\right)\!\!(w),\quad w\in \bndry D.
\end{equation*}
Then in fact
$$
d\l (w) \!=\!  \Lambda (w) d\sigma (w)
 $$
 with
 \begin{equation}\label{E:LLm}
\Lambda (w) = (n-1)!(4\pi)^{-n}|\det \r (w)|\ |\nabla\r (w)|,\quad w\in\bndry D,
\end{equation}
where $\det \r (w)$ is the determinant of the so-called ``Levi form for $D$'', which may be
 identified\footnote{here we have chosen a local coordinate system with respect to which the complex tangent space to $\bndry D$ at $w$ is identified with the space $\{(z_1,\ldots, z_{n-1}, 0)\, |\, z_j\in\C\}$.}
  with the 
 matrix
\begin{equation*}
\left\{\frac{\dee^2\r}{\dee z_j\dee\bar{z}_k}\right\}\!\bigg|_{z=w}\, ,\quad 1\leq j, k\leq n-1\, ,
\end{equation*}
see \cite[Lemma VII.3.9]{Ra-1}. Since $D$ is strongly pseudoconvex and of class $C^2$ it follows that
$$
c_1\leq \Lambda (w)\leq c_2\, ,\quad w\in\bndry D,
$$
so that $\l \approx \sigma$ (the two measures are mutually absolutely continuous) and thus the
$\Ctre$'s will be equivalently bounded with respect to either measure.

$\bullet$ Secondly, we have that the function
$$
\d (w, z) = |\langle\eta (w, z), w-z\rangle|^{1/2}\, ,\quad w, z\in\bndry D
$$
is a {\em quasi-distance}, namely, for any $w, z, \zeta \in \bndry D$ we have
$\d (w, z)= 0 \iff w=z$; $\d(w, z) \approx \d (z, w)$, and 
$\d (w, z) \leq C \big(\d(w, \zeta) + \d(\zeta, z)\big)$,  see \cite[Part I]{LS-5}.

$\bullet$  Moreover we have that the ensemble $X=\{\bndry D;\, \d; \l\}$ defines a {\em space of homogeneous type} with {\em homogeneous dimension $2n$}, see \cite[Part I]{LS-5}. That is, we have that $\l$ is a
doubling measure for the boundary balls $\B_r(w) =\{z\in\bndry D\, |\, \d(z, w)\leq r\}$ and in fact
$$
\l \big(\B_r(w)\big) \approx r^{2n}\, \quad \mbox{for any}\ w\in\bndry D\quad \mbox{and}\ r>0.
$$
\subsection{$L^2(\bndry D, d\l)$-regularity of the $\Ctre$'s: ad-hoc decompositions; the role of the Leray-Levi measure; application of $T(1)$.}\label{SS:2.8}
 In order to take full advantage of the measure $\l$ we make  the following decomposition of each of the $\Ctre$'s, see \cite[Part I]{LS-5}:
 $$
\Ctre = \Ctrse + \Rtre\, .
$$
Here the ``essential part'' $\Ctrse$ has kernel
$$
C_\epsilon^{\sharp}(w, z)\ =\ \frac{d\l (w)}{\langle\eta^{(\epsilon)}(w, z), w-z\rangle^n}\ =\
\frac{1}{(2\pi i)^n}\frac{j^*(\dee\rho\wedge(\deebar\dee\rho)^{n-1}(w)}
{\langle\eta^{(\epsilon)}(w, z), w-z\rangle^n}
$$
and captures the full singularity of $\Ctre$ in the sense that the ``remainders'' $\Rtre$'s  are smoothing operators that  map: $L^2(\bndry D)\to
C(\bar D)$, so in particular proving $L^2(\bndry D)$-regularity for $\Ctre$ is equivalent to proving the corresponding result for $\Ctrse$ (and we will henceforth ignore the $\Rtre$'s).

Now there is a further decomposition of $\Ctrse$, and a corresponding one for its formal adjoint on $L^2(\bndry D, d\l)$ that will play an important role in the application of the $T(1)$-theorem.
 The basic idea 
   is that one may express the kernels of each of $\Ctrse$ and its $L^2(\bndry D, d\l)$-adjoint 
   $(\Ctrse)^*$ as ``appropriate derivatives'' (plus acceptable remainders) that is,  as the differentials of $(2n-2)$-forms whose coefficients have better homogeneity than 
the kernels of each of $\Ctrse$ and $(\Ctrse)^*$; the desired decompositions will then result by an application of Stokes' theorem.  To put matters more precisely, given $f\in C^1(\bndry D)$ we have \cite[Part I]{LS-5}
$$
\Ctrse (f) (z) = \Etre(df)(z) +\Rtrse(f)(z),\quad \mbox{for}\ z\in \bndry D\, 
$$
and
$$
(\Ctrse)^{\!*} (f) (z) = \widetilde\Etre(df)(z) +\widetilde\Rtrse(f)(z),\quad \mbox{for}\ z\in \bndry D\, .
$$
Here the ``essential parts'' $\Etre$ and $\widetilde\Etre$ act on continuous 1-forms 
$\omega$ on $\bndry D$
as follows:
    $$
    \Etre(\omega)(z) = 
   c_n \int_{w\in\bndry D}
      \frac{\omega\wedge j^*(\deebar\dee\rho)^{n-1}(w)}
   {\langle\eta^{(\epsilon)}(w, z), w-z\rangle^{n-1}}, \quad z\in\bndry D,
    $$
with $c_n =1/[(n-1)(2\pi i)^n]$. Comparing $\Etre$ with $\Ctrse$ we see that the
kernel of $\Etre$ has the
 improved homogeneity 
${\langle\eta^{(\epsilon)}(w, z), w-z\rangle^{-n+1}}$ (as opposed to 
${\langle\eta^{(\epsilon)}(w, z), w-z\rangle^{-n}}$).
Similarly, we have
    $$
    \widetilde\Etre(\omega)(z) = 
   c_n \int_{w\in\bndry D}
      \frac{\omega\wedge j^*(\deebar\dee\rho)^{n-1}(w)}
   {\langle\bar\eta^{(\epsilon)}(z, w), \bar z-\bar w\rangle^{\,n-1}}, \quad z\in\bndry D.
    $$
In both decompositions, the remainders $\Rtrse$ and $\widetilde\Rtrse$ are once again
smoothing operators: $L^2(\bndry D, d\l)\to C(\bar D)$.
We point out that the decomposition for $(\Ctrse)^*$ takes full advantage
of the measure $d\l$, in the sense that the decomposition for the adjoint of $\Ctrse$ is valid only if the 
adjoint
 is computed with respect to the duality for $L^2(\bndry D, d\l)$ (there is no such decomposition for the adjoint of $\Ctrse$ in $L^2(\bndry D, d\sigma))$. 

Using these decompositions one then shows that the functions
\begin{equation}\label{E:holder}
\h = \Ctrse (1)\quad \mbox{and}\quad \h^* = (\Ctrse)^* (1)
\end{equation}
are continuous on $\bndry D$, see \cite[Part I]{LS-5}.

The properties of the homogenous space $\{\bndry D, \d, \l\}$ along with the above decompositions for $\Ctrse$ and its $L^2(\bndry D, d\l)$-adjoint $(\Ctrse)^*$ now ensure that for any $0<\epsilon\leq \epsilon_0$, the operators
$T:= \Ctrse$ satisfy all the hypotheses of the $T(1)$-theorem, namely: the kernel of $\Ctrse$ satisfies 
appropriate size and regularity estimates; the operator $\Ctrse$ is weakly bounded and satisfies the cancellation conditions\footnote{In fact $\h$ and $\ha$ can be proved to be be H\"older continuous in the sense that e.g., $|\h (w) -\h (z)|\leq C\d(w, z)^\alpha,\ w, z\in\bndry D$ for any $0<\alpha<1$, and this in turn allows to reduce the application of the $T(1)$ theorem to the simpler setting: $T_0(1)=0=T^*_0(1)$ for a suitable auxiliary operator $T_0$, see \cite[Section 6.3]{LS-4} and also \cite[Part I]{LS-5}.} \eqref{E:holder}, see \cite[Part I]{LS-5}. 
 This concludes the proof of the $L^p(\bndry D, d\l)$-regularity of $\Ctrse$ 
 and therefore the proof of item {\tt (c.)} (see Section
 \ref{SS:2.3})
  for each of the operators $\Ctre$'s.
 
 \subsection{$L^p(\bndry D, d\l)$-regularity of the Szeg\H o projection (Leray-Levi measure)}\label{SS:2.9} The notion of orthogonal projection (in particular the definition of the Szeg\H o projection)
 relies on the specific measure 
  that is being used in the definition of  $L^2(\bndry D)$: different measures give rise to different orthogonal projections and there no simple way of deducing regularity for one projection from the corresponding result for the other.  In this section we highlight the procedure carried out in \cite[Part II]{LS-5} to solve the $L^p(\bndry D)$-regularity problem for the Szeg\H o projection defined with the respect to the Leray-Levi measure and its corresponding duality on $L^2(\bndry D)$:
 \begin{equation}\label{E:gen-dual}
 (f, g) =\int\limits_{w\in\bndry D}\!\!\!f(w)\bar g(w)\, d\l (w).
 \end{equation}
We will denote such projection $\Sl$. What is still missing from the procedure summarized in Section \ref{SS:2.3} is item {\tt (d.)}, namely the ``smallness'' of the operators $\Aope =
\Ctre^*-\Ctre$, where the adjoint $\Ctre^*$ is computed with respect to the duality
\eqref{E:gen-dual}.
 Going back to  the setting of \cite{KS-1} ($D$ of class $C^3$) such smallness
resulted from the following estimate for the kernel of $\mathcal A^{(\sigma)}$ (denoted below by 
$A^{(\sigma)} (w, z)$)
\begin{equation}\label{E:est-6}
|A^{(\sigma)} (w, z)|\leq C|w -z|^{2+\beta} \quad \mbox{whenever}\quad \d(w, z)\leq \delta
\end{equation}
for some $\beta>0$ (in fact for $\beta =1$), which ultimately gave the compactness of $\mathcal A^{(\sigma)}$ in $L^p(\bndry D, d\sigma)$: one considered the 
 operators $\{\mathcal A^{(\sigma)}_\delta\}_{_{\!\delta}}$ with kernels 
$$\left( 1-\chi_\delta(|w-z|^2) \right)A^{(\sigma)} (w, z)$$
for a smooth cutoff function $\chi_\delta(t)$. Such operators are obviously compact in $L^p(\bndry D, d\sigma)$ for $1<p<\infty$. Now the estimate \eqref{E:est-6} grants 
$\| \mathcal A^{(\sigma)}_\delta -
\mathcal A^{(\sigma)}\|_{L^p\to L^p}\leq C \delta^\beta$, and the compactness of  $\mathcal A^{(\sigma)}$ then follows by letting $\delta \to 0$.

  Note that the positivity of $\beta$ in \eqref{E:est-6} is crucial:
   the estimate
 $|A^{(\sigma)} (w, z)|\leq C\d(w, z)^2$
would only yield the inconclusive inequality $\| \mathcal A^{(\sigma)}_\delta -
\mathcal A^{(\sigma)}\|_{L^p\to L^p}\leq C$.
However in our less regular setting there is no analog of \eqref{E:est-6} with $\beta>0$. In fact
the operator $\Aope$ will in general fail to be compact in $L^p(\bndry D, \lambda)$, see 
\cite[Corollary 5]{BaLa},
 and one must proceed by a different analysis. 
 What holds in place of \eqref{E:est-6} is the following, ``weaker'' smallness for the kernels  of the operators $\Aope$ (denoted $A^{(\lambda)}_\epsilon (w, z)$), namely:
 \begin{equation}\label{E:est-7}
|A^{(\lambda)}_\epsilon (w, z)|\leq C\,\epsilon\, \d(w, z)^2 \quad \mbox{whenever}\quad \d(w, z)\leq \delta\, ,
\end{equation}
To use \eqref{E:est-7} we consider the operators 
$\mathcal A^{(\lambda)}_{\epsilon, \delta}$ with kernel
$$
\chi_\delta(\d(w,z))\, A^{(\lambda)}_\epsilon (w, z)
$$
where $A^{(\lambda)}_\epsilon (w, z)$ is the kernel of $\Aope$ and $\chi_\delta(t)$ is a smooth cutoff function.
Then we have that
\begin{equation}\label{E:decA}
\Aope = \mathcal A^{(\lambda)}_{\epsilon, \delta} + \mathcal R^{(\lambda)}_{\epsilon, \delta}\, .
\end{equation}
Now one may apply the $T(1)$-theorem (as in the previous section) to prove that 
the ``essential part'' $\mathcal A^{(\lambda)}_{\epsilon, \delta}$ is bounded: 
$L^p(\bndry D, d\l) \to L^p(\bndry D, d\l)$ for $1<p<\infty$. In fact \eqref{E:est-7}
yields the improved estimate \cite[Part II]{LS-5}
\begin{equation}\label{E:est-8}
\|\mathcal A^{(\lambda)}_{\epsilon, \delta}\|_{L^p(\bndry D, d\l) \to L^p(\bndry D, d\l)}\leq 
\epsilon^{1/2}M_p,\quad 1<p<\infty
\end{equation}
for any $0<\delta<\delta_0 (\epsilon)$ and for any $0<\epsilon\leq \epsilon_0$, where the bound $M_p$ is symmetric in $p$, i.e. $M_p=M_{p'}$ whenever $1/p + 1/p'=1$.

On the other hand, the ``remainder'' operators
 $\mathcal R^{(\lambda)}_{\epsilon, \delta}$ (whose kernel are supported outside of the critical balls $\{\d(w, z)<\delta\}$)
are readily seen to map: $L^1(\bndry D, d\l)\to L^\infty(\bndry D)$ (although their $L^p\to L^p$-norms may be very large). 

We now proceed to compare $\Sl$ with the Cauchy-type integrals $\Ctre$.
By items {\tt (a.)} and {\tt (b.)} (proved in Section \ref{SS:2.6}) and proceeding as in Section \ref{SS:2.3}, we recover the identity
\begin{equation}\label{E:orig-id}
\Sl\, [I-\Aope]=\Ctre\quad \mbox{in}\ L^2(\bndry D, d\l)\quad \mbox{for any}\ \ 0<\epsilon<\epsilon_0.
\end{equation}
 Combining the above with \eqref{E:decA} we get
 $$
 \Sl\, [I-\mathcal A^{(\lambda)}_{\epsilon, \delta} ]= \Ctre - \Sl\,\mathcal R^{(\lambda)}_{\epsilon, \delta}\quad 
 \mbox{in}\ L^2(\bndry D, d\l)\quad \mbox{for any}\ \ 0<\epsilon<\epsilon_0.
 $$
 We now fix $1<p<\infty$ and prove  $L^p(\bndry D, d\lambda)$-regularity of $\Sl$ for such $p$; for the time being we take $1<p<2$, so that the two inclusions: $L^p(\bndry D, d\l)\hookrightarrow L^1(\bndry D, d\l)$ and $L^2(\bndry D, d\l)\hookrightarrow L^p(\bndry D, d\l)$ are bounded\footnote{here we are using the hypothesis that the domain $D$ is bounded.}. It follows that the operator
 $$
 \Sl\,\mathcal R^{(\lambda)}_{\epsilon, \delta}: L^p\hookrightarrow L^1\to
  L^\infty\hookrightarrow L^2\to L^2\hookrightarrow L^p
 $$
 is bounded for any $0<\delta<\delta_0(\epsilon)$ and any $0<\epsilon\leq \epsilon_0$ (here we have used the $L^1\to L^\infty$-regularity of $\mathcal R^{(\lambda)}_{\epsilon, \delta}$ and the $L^2\to L^2$-regularity of $\Sl$ ). Moreover
 $\Ctre$ is bounded: $L^p\to L^p$ for any $0<\epsilon\leq \epsilon_0$ by item {\tt (c.)} (which was proved in sections \ref{SS:2.7} and 
 \ref{SS:2.8}). We now fix $\epsilon = \epsilon (p)<<1$ so that
 $$
 \epsilon^{1/2} M_p<1\, ,
 $$
where $M_p$ is as in \eqref{E:est-8}. Then by \eqref{E:est-8}  we have that for any $\delta\leq \delta_0(\epsilon)$ the operator
 $$
I- \mathcal A^{(\lambda)}_{\epsilon, \delta}
 $$
 is invertible in $L^p(\bndry D, d\l)$ by a partial Neumann series, and has bounded inverse.
 
 We conclude from the above that
 \begin{equation}\label{E:KS-2}
 \Sl\,=\, [\Ctre - \Sl\,\mathcal R^{(\lambda)}]\,[I-\mathcal A^{(\lambda)}_{\epsilon, \delta} ]^{-1}\quad \mbox{in}\quad L^2(\bndry D, d\l).
 \end{equation}
 However, by what has been said, the right-hand side of this identity is a bounded operator in $L^p(\bndry D, d\l)$, thus showing that $\Sl$ extends to a bounded operator in 
 $L^p(\bndry D, d\l)$ for any $1<p\leq 2$. By duality (and the fact that $(\Sl)^* = \Sl$)  it follows that $\Sl$ is also bounded in $L^p(\bndry D, d\l)$ for any $2\leq  p<\infty$.
 The $L^p$-regularity problem for $\Sl$ is therefore solved with $P=\infty$, whenever $D$ is a bounded, strongly pseudoconvex domain of class $C^2$.
 \subsection{$L^p(\bndry D, d\l)$-regularity of the Szeg\H o projection: other measures}\label{SS:2.10}
 We recall from Section \ref{SS:2.7} that the Leray-Levi measure $\l$ and the induced Lebesgue
  measure $\sigma$ are mutually absolutely continuous,
    see \eqref{E:LLm} and comments thereafter. It follows that the  holomorphic Cauchy-type integrals 
  $\{\Ctre\}_\epsilon$ are equivalently bounded in $L^p(\bndry D, d\l)$ and $L^p(\bndry D, d\sigma)$ and, more generally in  $L^p(\bndry D, \varphi d\l)$ where $\varphi$ is any 
  continuous function on $\bndry D$ with uniform upper and lower bounds. On the other
  hand, if we denote the Szeg\H o projection for $L^2(\bndry D, \varphi d\l)$ by $\Sp$, there is no direct way to compare $\Sl$ with $\Sp$ that would allow to deduce $L^p(\bndry D, \varphi d\l)$-regularity for $\Sp$ from $L^p(\bndry D, d\l)$-regularity of $\Sl$.

  Instead, what we are able to compare the two adjoints of the operator $\Ctre$ with respect
  to these different measures. Letting $\mathcal A^{(\varphi)}_{\epsilon, \delta}$ denote
  the operator analogous to \eqref{E:decA} (however now defined for the measure $\varphi d\l$), we have the following inequality \cite[Part II]{LS-5}
\smallskip
  
  $\bullet$
  \quad 
  $\displaystyle{\|\mathcal A^{(\varphi)}_{\epsilon, \delta}\|_{L^p\to L^p} \leq
  \|\mathcal A^{(\lambda)}_{\epsilon, \delta}\|_{L^p\to L^p}
  +
  \sup\limits_{w\in\bndry D}|\varphi^{-1}(w)|\, \|[\mathcal C_{\epsilon, \delta}, \varphi]\|_{L^p\to L^p}}$
  \medskip
  
  where $\mathcal C_{\epsilon, \delta}$ is the operator with kernel $\chi_\delta(\d(w, z))\, C^{(\epsilon)}(w, z)$, see \eqref{E:LS-ker}, and 
  $[T, \varphi] = T\,M_\varphi - M_\varphi T$ denotes the commutator with the multiplication operator $M_\varphi$ (multiplication by $\varphi$). 
  
  Furthermore, we have that \cite[Part II]{LS-5}
  
  $\bullet$
  \quad 
  $\displaystyle{\| [\mathcal C_{\epsilon, \delta}, \varphi]\|_{L^p\to L^p}}\leq \epsilon M_p$
for any $\delta \leq \delta_0(\epsilon)$ and for any $0<\epsilon\leq \epsilon_0$.
 \medskip
 
 Taking these two facts into account, the proof of the $L^p(\bndry D, \varphi d\l)$-regularity of $\Sp$ may now be obtained by following the same steps as in the proof of the corresponding result for $\Sl$ on $L^p(\bndry D, d\l)$, see Section \ref{SS:2.9}.
 
 \section{Further results}
\subsection{The Bergman projection}\label{SS:3.1} One may also state the $L^p$-regularity problem
for the {\em Bergman projection}, that is the orthogonal projection of $L^2(D, dV)$ onto the Bergman space $\vartheta (D)\cap L^2(D, dV)$ (namely, the space of functions that are holomorphic in $D$ and square-integrable on $D$ with respect to the measure on $D$ induced by the Lebesgue measure for $\Bbb \C^n$ via the inclusion: $D\subset \Cn$).
The $L^p$-regularity problem for the Bergman projection was studied by Ligocka \cite{Li} in the setting of bounded, strongly pseudconvex domains of class $C^4$, and was recently extended \cite{LS-2} to the class $C^2$, with $L^p(D, dV)$-regularity holding for $1<p<\infty$. This problem can be approached in a fashion similar to the $L^p$-regularity problem for the Szeg\H o projection, but is in fact simpler than that problem, in several respects:

$\bullet$ There is no advantage in considering ``ad-hoc''  volume measures for the domain $D$ (some ``solid'' version of the Leray-Levi measure) and one may work directly with the induced Lebesgue measure $dV$. 

$\bullet$ In this context, the role of the ``holomorphic Cauchy integrals'' $\Ctre$'s is played by
``solid'' integral operators $\bf B_{\epsilon}$ 
 acting on $L^p(D, dV)$, whose kernel
is essentially the ``derivative'' of the kernels of the $\Ctre$'s, specifically, it is the
$2n$-form $\deebar_w \widetilde C^{(\epsilon)}(w, z)$, see \eqref{LS-kere} (then corrected to achieve global holomorphicity). Such operators will still produce (and reproduce) holomorphic functions from merely $L^1$ data, see \cite[Propositions 3.1 and 3.2]{LS-3}. By their nature, these operators are less singular than the $\Ctre$'s and their $L^p (D, dV)$-regularity can be established by 
direct means (with no need to invoke the $T(1)$-theorem, \cite[Section 4]{LS-2}).

$\bullet$ The solution of the $L^p$-regularity problem for the Bergman projection of a strongly pseudoconvex domain of class $C^2$ 
now follows a parallel argument to the corresponding result for the Szeg\H o projection,
by proving ``$\epsilon$-smallness'' for the kernels of $\bf B_{\epsilon}^*- \bf B_{\epsilon}$
with $\epsilon$ again tailored to the size of the Lebesgue exponent $p$, see \cite[Sections 5. and 6.]{LS-2}.

$\bullet$ In fact one also proves $L^p$-regularity for the ``absolute Bergman projection'',
that is the operator whose kernel is the {\em absolute value} of the Bergman kernel, see \cite[Section 6]{LS-2} and \cite[Section I.1, Example 1]{C-2}. (We point out that the corresponding statement for the ``absolute Szeg\H o projection'' is known to be false by the very nature of the Szeg\H o kernel, whose treatment requires cancellation conditions that would be lost by considering its absolute value.)
\subsection{Holomorphic Cauchy integrals below the $C^2$-threshold}\label{SS:3.2}
 A  theory of  holomorphic Cauchy integrals can also be developed for so-called  {\em strongly $\C$-linearly convex} domains.
 While $\C$-linear convexity is a stronger notion than pseudoconvexity (in the sense that any strongly $\C$-linearly convex domain
of class $C^2$ is strongly pseudoconvex  but the converse is not true \cite[ Proposition 3.2 and Example pg.797]{LS-4}), it is a notion that rests on only
{\em one} derivative of the defining function and is therefore naturally meaningful for domains of class $C^1$. See \cite{APS}, \cite{Ho} and \cite[Section 3]{LS-4} for the definition and main properties of $\C$-linear convexity.

In \cite{LS-4} we study existence and regularity of holomorphic Cauchy integrals for bounded, strongly $\C$-linear convex domains of class $C^{1,1}$.
In this context the $C^{1,1}$ category plays a role analogous
to the Lipschitz category for a planar domain; the relevant kernel is the {\em Cauchy-Leray kernel}
\begin{equation}\label{E:CLker}
K(w, z) = \frac{1}{(2\pi)^n}j^*\left(\!\frac{\dee\rho\wedge (\deebar\dee\rho)^{n-1}(w)}
{\langle\dee\rho (w), w-z\rangle^n}\!\right) = \frac{d\l (w)}{\langle\dee\rho (w), w-z\rangle^n}\, .
\end{equation}

$\bullet$ This kernel was first identified by Leray \cite{L} in the context of (strongly convex) domains of class $C^2$. A new and substantial obstacle that arises in the $C^{1,1}$
  setting is the fact that the (familiar)
numerator of $K(w, z)$ may not make sense:
 while the Rademacher Theorem ensures
that the $C^{1,1}$ function $\rho$ be twice differentiable almost everywhere in 
$\Cn$, here we are taking its restriction to the boundary $\bndry D$ which is in fact 
a zero-measure subset of $\Cn$, and the coefficients of $\deebar\dee\rho$ may indeed be undefined on $\bndry D$ (explicit examples can be given); however, the pullback by the inclusion $j^*(\deebar\dee\rho)$ only pertains the {\em tangential components} of such coefficients, which are indeed well-defined, see \cite[Proposition 23]{LS-4}. As a result, one has that the Leray-Levi measure $d\l$ is well defined also in this less regular context, and it is again equivalent to the induced Lebesgue measure $d\sigma$, see \cite[Section 3.4]{LS-4}.

$\bullet$ The strong $\C$-linear convexity of $D$ ensures that the global bound: 
$|\langle\dee\rho (w), w-z\rangle|\geq c|w-z|^2$ holds for any $w\in \bndry D$ and any $z\in \bar D$, see \cite[(3.4)]{LS-4}. Thus, the Cauchy-Leray kernel is globally holomorphic (no need for correction).
While not universal in the sense of Section \ref{SS:2.4}, The Cauchy-Leray kernel is {\em canonical} in the sense that it does not depend on the choice of 
defining function (while each of the numerator and denominator in \eqref{E:CLker} depend on the choice of defining function $\rho$, their quotient does not, see \cite[Proposition 4.1]{LS-4}). This is  in great contrast
with the situation for the kernels of the operators $\Ctr$ and $\Ctre$ considered in sections \ref{SS:2.5} -- \ref{SS:2.8},
  which do depend on the choice of defining function and are thus non-canonical.

$\bullet$ Letting $\mathcal K$ denote the Cauchy-Leray operator with kernel \eqref{E:CLker}, we prove that  $\mathcal K$ is bounded in $L^p(\bndry D, d\l)$  for any $1<p<\infty$ (and thus $L^p(\bndry D, d\sigma)$), by a $T(1)$-theorem for a space of homogeneous type again informed 
by the geometry and regularity of $\bndry D$ (in a spirit that is similar to the situation described in Sections \ref{SS:2.7} and \ref{SS:2.8}), and again under the simpler cancellation conditions $T(1)=0= T^*(1)$, see \cite[Section 6]{LS-4}. 

$\bullet$  Our hypotheses on the domain $D$ are optimal in the sense that for any $0<\alpha<1$ there are strongly $\C$-linearly convex domains 
 $D_\alpha$ of class $C^{1, \alpha}$, for which the Cauchy-Leray operator $\mathcal K$ is well-defined but unbounded on each of $L^2(\bndry D, d\l)$ and $L^2(\bndry D, d\sigma)$, see \cite[Section 6, Example 2]{BaLa}; similarly, there is a smooth {\em weakly} $\C$-linearly convex domain $D$ for which $\mathcal K$ is well-defined but unbounded on each of $L^2(\bndry D, d\l)$ and $L^2(\bndry D, d\sigma)$, see 
 \cite[Section 6, Example 1]{BaLa}.
  
$\bullet$ The difference $\mathcal K^* - \mathcal K$ has no inherent ``smallness'':
not even the weaker ``$\epsilon$-smallness'' \eqref{E:est-7} (if, say, one were to approximate $\rho$ with smoother functions $\tau^{(\epsilon)}$)\footnote{This failure  occurs essentially because the denominator 
$\langle\dee\rho (w), w-z\rangle^n$ is ``unrefined'' in the sense that it uses only one derivative of the defining function.}. Thus the study of the Bergman and Szeg\H o projections for strongly $\C$-linearly convex domains requires a different approach
and is the object of current investigation.

\subsection{Representations for the Hardy and Bergman spaces of holomorphic functions} As an application of the $L^p$-regularity of the holomorphic Cauchy-type integrals and of the Szeg\H o and Bergman projections we obtain various representations for the Hardy and Bergman
spaces of holomorphic functions. Specifically, for a strongly pseudoconvex domain of class $C^2$, we have the following, see \cite{LS-6} (see also \cite{La-2}):

$\bullet$ The space of functions holomorphic in a neighborhood of $\bar D$ is dense in 
$\mathcal H^p(\bndry D, \varphi d\l)$ (a consequence of the $L^p$-regularity of the $\Ctre$'s).

$\bullet$ $\Ctre: L^p(\bndry D, \varphi d\l) \to \mathcal H^p(\bndry D, \varphi d\l)$ for any
$1<p<\infty$. Furthermore, $f\in  \mathcal H^p(\bndry D, \varphi d\l)$ if, and only if
$\Ctre f = f$ (again a consequence of the $L^p$-regularity of the $\Ctre$'s).
 
 Corresponding statements hold for the situation when $D$ is strongly $\C$-linearly convex and of class $C^{1,1}$ (with the $\Ctre$'s replaced by the Cauchy-Leray operator $\mathcal K$).

Furthermore, for a strongly pseudoconvex domain of class $C^2$ we have, see \cite[Proposition 7.1]{LS-3} and \cite{LS-6}:

$\bullet$ The space of functions holomorphic in a neighborhood of $\bar D$ is dense in 
$\vartheta (D)\cap  L^p(D, dV)$ (a consequence of the $L^p$-regularity of the $\mathbf B_{\epsilon}$'s).

$\bullet$ $\Sp: L^p(\bndry D, \varphi d\l) \to \mathcal H^p(\bndry D, \varphi d\l)$ for any
$1<p<\infty$. Moreover, $f\in  \mathcal H^p(\bndry D, \varphi d\l)$ if, and only if
$\Sp f = f$ (a consequence of the $L^p$-regularity of $\Sp$).

 \bibliographystyle{alpha}

\end{document}